\documentclass[11pt]{amsart}
\usepackage{amsmath}
\usepackage{amssymb}
\usepackage{epsfig}
\usepackage{latexsym}
\usepackage{amscd}
\newtheorem{thm}{Th\'eor\`eme}
\newtheorem{cor}{Corollaire}
\newtheorem{prop}{Proposition}
\newtheorem{defi}{D{\'e}finition}

\newtheorem{ex}{Exemple}
\title{Big{\`e}bres quasi-Lie et boucles de Lie}
\author{Momo BANGOURA}
\address{D{\'e}partement de Math{\'e}matiques, Universit{\'e} de Conakry,\\
BP 1147, R{\'e}publique de Guin{\'e}e}
\email{angoura@gn.refer.org}
\thanks{Visiting Scientist at The Abdus Salam ICTP\\
{2000 {\it Mathematics Subject Classification}. 17B70
17A30.
{\it Mots cl{\'e}s. Action, alg{\`e}bre de Lie, alg{\`e}bre
d'Akivis, big{\`e}bre quasi-Lie, boucle de Lie mono-alternative
{\`a} droite, groupe de Lie, quasi-big{\`e}bre de Lie,
quasi-double, quasi-groupe, quasi-poisson}}}
\date{}
\begin{document}
\begin{abstract}
In this work, we define the quasi-Poisson  Lie quasigroups, dual
objects to the quasi-Poisson Lie groups and we establish the
correspondance  between the local quasi-Poisson Lie quasigoups and
quasi-Lie bialgebras (up to isomorphism).

\vspace{0.5cm} {\bf R{\'e}sum{\'e} :} Dans ce travail, nous
d{\'e}finissons les quasi-groupes de Lie quasi-Poisson, objets
duaux des groupes de Lie quasi-Poisson et nous {\'e}tablissons une
correspondance biunivoque ({\`a} isomorphisme pr{\`e}s) entre les
quasi-groupes de Lie quasi-Poisson locaux et les big{\`e}bres
quasi-Lie.
\end{abstract}
\maketitle
\tableofcontents
\section{Introduction}
Les groupes de Lie-Poisson {\'e}tant les limites classiques des
alg{\`e}bres de Hopf, les groupes de Lie quasi-Poisson \cite{KS1,
KS2} sont les limites classiques des quasi-alg{\`e}bres de Hopf
introduites par Drinfeld \cite{Dr2}, et les objets
infinit{\'e}simaux sont appel{\'e}s les quasi-big{\`e}bres de Lie
dans \cite{Dr2} et quasi-big{\`e}bres jacobiennes dans \cite
{KS2}. Contrairement aux notions de big{\`e}bre de Lie et groupe
de Lie-Poisson \cite{Dr1, KS-Ma, A-KS, Lu-W, Lu}, les notions de
quasi-big{\`e}bre de Lie et de groupe de Lie quasi-Poisson ne sont
pas auto-duales; l'objet dual d'une quasi-big{\`e}bre de Lie est
appel{\'e} une quasi-big{\`e}bre co-jacobienne dans \cite{KS2,
Ba-KS}. Nous l'appellerons big{\`e}bre quasi-Lie. Ici, nous nous
proposons de construire l'objet dual d'un groupe de Lie
quasi-Poisson, c'est-{\`a}-dire un objet g{\'e}om{\'e}trique
g{\'e}n{\'e}ralisant les groupes de Lie-Poisson dont l'espace
tangent en un {\'e}l{\'e}ment distingu{\'e} est muni d'une
structure de big{\`e}bre quasi-Lie et nous convenons d'appeler cet
objet quasi-groupe de Lie quasi-Poisson. Plus pr{\'e}cisement, un
quasi-groupe de Lie quasi-Poisson est un $G-$espace quasi-Poisson
au sens de \cite{Al-KS}, muni d'une structure de boucle de Lie
mono-alternative {\`a} droite avec certaines relations de
compatibilit{\'e} entre les deux structures, notamment celle
g{\'e}n{\'e}ralisant la multiplicativit{\'e}(\cite{Lu-W, Lu}) du
champ de bivecteurs qui est satisfaite dans le cas des groupes de
Lie-Poisson. Pour cela, on se sert de la th{\'e}orie des boucles
de Lie mono-alternatives {\`a} droite \cite{Sa-Mi2} en
introduisant les notions de groupe de Lie quasi-double et
d'alg{\`e}bre de Lie quasi-double dont les propri{\'e}t{\'e}s font
ressortir l'action d'un groupe de Lie sur une boucle de Lie et
d'une action par d{\'e}faut de cette boucle de Lie sur le groupe
de Lie. Un groupe de Lie quasi-double est un triplet form{\'e}
d'un groupe de Lie, d'un sous-groupe de Lie ferm{\'e} et d'une
sous-vari{\'e}t{\'e} ferm{\'e}e du groupe de Lie, laquelle est
munie d'une structure de boucle de Lie mono-alternative {\`a}
droite, le triplet v{\'e}rifiant un certain nombre de
propri{\'e}t{\'e}s; l'objet infinit{\'e}simal associ{\'e} est
app{\'e}l{\'e} une alg{\`e}bre de Lie quasi-double. Cette
construction g{\'e}n{\'e}ralise certains r{\'e}sultats de
\cite{Lu-W} sur les groupes de Lie doubles et alg{\`e}bres de Lie
doubles; ces derni{\`e}res sont aussi appel{\'e}es alg{\`e}bres de
Lie bicrois{\'e}es dans \cite{KS-Ma}.\\
Dans la section 2 nous faisons un rappel de quelques notions sur
les boucles de Lie et les alg{\`e}bres d'Akivis.\\
Dans la section 3, on d{\'e}finit et on {\'e}tudie les groupes de
Lie quasi-doubles et les alg{\`e}bres de Lie quasi-doubles en
{\'e}tablissant une correspondance entre les deux notions tout en
faisant ressortir la compatiblit{\'e} entre les diff{\'e}rentes
op{\'e}rations les d{\'e}finissant. \\
Dans la section 4, on fait un bref rappel sur les big{\`e}bres
quasi-Lie en montrant leur lien d'une part avec les alg{\`e}bres
de Lie quasi-doubles, et d'autre part avec les alg{\`e}bres d'Akivis.\\
La derni{\`e}re section recouvre l'{\'e}ssentiel du travail {\`a}
savoir la d{\'e}finition des quasi-groupes de Lie quasi-Poisson et
l'{\'e}tablissement d'une correspondance bijective entre les
quasi-groupes de Lie quasi-Poisson locaux et les big{\`e}bres
quasi-Lie.
\section{Boucles de Lie mono-alternatives {\`a} droite}
Dans cette section, nous faisons un rappel sur les boucles de Lie
en reprenant les d{\'e}finitions et certains r{\'e}sultats de
\cite{Sa, Sa-Mi1, Sa-Mi2}.
\begin{defi} Une boucle de Lie est une vari{\'e}t{\'e} analytique
r{\'e}elle $B$ munie d'une op{\'e}ration interne analytique
$$\begin{array}{ccccc}
m &:& B\times B& \longrightarrow & B \\
  & & (a, b)& \longrightarrow & $a.b = m(a, b)$
\end{array}$$
admettant un {\'e}l{\'e}ment neutre $\varepsilon$,
v{\'e}rifiant
$$a.\varepsilon = \varepsilon .a = a$$
pour tout $a\in B$, telle que les {\'e}quations
\begin{center}$a.b
= c$ \hspace{0.5cm} et \hspace{0.5cm}$b.a = c$\end{center}
admettent des solutions uniques (not{\'e}es respectivement
$a\setminus c$ et $c/a$) et les applications
\begin{center}$(a, b)\longrightarrow b\setminus a$, \hspace{1cm}$(a, b)\longrightarrow
a/b$\end{center} sont analytiques.
\end{defi}
Ainsi, comme on le voit, tout groupe de Lie est une boucle de
Lie.\\
Si les trois op{\'e}rations \\
\begin{center}$(a, b)\longrightarrow ab$, \hspace{0.5cm} $(a,
b)\longrightarrow a/b$ , \hspace{0.5cm} $(a, b) \longrightarrow
b\setminus a$ \end{center}
sont d{\'e}finies localement au
voisinage de $\varepsilon$, on dit
que $B$ est une boucle de Lie locale.\\
{\bf Remarque :} L'unicit{\'e} des solutions dans une boucle de
Lie des {\'e}quations $a.b = c$ et $b.a = c$ signifie que les
translations {\`a} gauche et {\`a} droite par un {\'e}l{\'e}ment
$a\in B$, not{\'e}es respectivement par
$\lambda_{a}$ et $\rho_{a}$, sont des diff{\'e}omorphismes.\\
L'espace tangent en l'identit{\'e} {\`a} toute boucle de Lie est
muni d'une structure g{\'e}n{\'e}ralisant la notion de structure
d'alg{\`e}bre de Lie, et dont les propri{\'e}t{\'e}s d{\'e}coulent
de celles de la boucle de Lie. Ce qui nous conduit {\`a} la
d{\'e}finition suivante :
\begin{defi} Une alg{\`e}bre d'Akivis $(A, [.,.], <.,.,.>)$ est
un espace vectoriel $A$ muni d'une application bilin{\'e}aire
antisym{\'e}trique\\
$$\begin{array}{ccccc}
[.,.]& : & A\times A & \longrightarrow & A\\
     &   & (x, y)& \longrightarrow & [x, y]
\end{array}$$
et d'une application trilin{\'e}aire\\
$$\begin{array}{ccccc}
<.,.,.> & : & A\times A\times A & \longrightarrow & A\\
     &   & (x, y, z)& \longrightarrow & <x, y, z>
\end{array}$$
telles que
$$\sum_{\sigma \in S_3}(Sign\sigma)<x_{\sigma(1)}, x_{\sigma(2)},
x_{\sigma_(3)}> = \sum_{\bigcirc}[[x_1, x_2], x_3]$$ o{\`u} $S_3$
d{\'e}signe le groupe des permutations de $(1, 2, 3)$,
$Sign\sigma$ la signature de $\sigma \in S_3$ et $\sum_{\bigcirc}$
d{\'e}signe la somme sur les permutations circulaires de $x_1,
x_2, x_3$.
\end{defi}
Ainsi, toute alg{\`e}bre de Lie d{\'e}termine une structure
d'alg{\`e}bre d'Akivis en prenant $<.,.,.> = 0$. \\
Pour une boucle de Lie locale $B$, on montre dans \cite{Du} qu'il
existe une structure d'alg{\`e}bre d'Akivis sur
$T_{\varepsilon}B$, l'espace tangent en l'identit{\'e},
d{\'e}termin{\'e}e par les op{\'erations} suivantes:
$$[x, y] =
\frac{1}{2}\frac{d^{2}}{dt^{2}}((x(t)y(t))/((y(t)x(t))|_{t=0}$$
$$<x, y, z> =
\frac{1}{6}\frac{d^{3}}{dt^{3}}(((x(t)y(t))z(t)/(x(t)(y(t)z(t)))|_{t=0}$$
pour $x, y, z \in T_{\varepsilon}B$, o{\`u} $x(t), y(t), z(t)$
sont des courbes de classe $C^{k}$ $(k\geq 3)$ dans $B$, passant
par le point $\varepsilon$ avec les vecteurs tangents $x, y, z$
respectivement.\\
L'alg{\`e}bre d'Akivis ainsi d{\'e}crite est appel{\'e}e
alg{\`e}bre d'Akivis de la boucle de Lie $B$.\\
Nous allons {\`a} pr{\'e}sent d{\'e}finir une classe de boucles de
Lie locales avec laquelle nous travaillerons dans la suite.
\begin{defi} Une boucle de Lie locale $B$ est dite mono-alternative {\`a}
droite, si pour tous entiers $k,l\in\mathbb{Z}$ et pour tous
$a,b\in B$, voisins de $\varepsilon$, on a :
$$(ab^{k})b^{l} = ab^{k+l}.$$
\end{defi}
{\bf Remarque :} La formule ci-dessus exprime d'une part une
associativit{\'e} faible de la multiplication d{\'e}finie sur $B$,
d'autre part l'existence d'un inverse pour tout {\'e}l{\'e}ment
$a\in B$, not{\'e} $a^{-1}$. Par convention , on a pour tout $a
\in B$, $a^{0} = \varepsilon$, $a^{-k} = (a^{-1})^{k}$ et donc
$\rho_{a}^{-1} = \rho_{a^{-1}}$.
\begin{ex} Soit $SH(n)$ l'ensemble des matrices hermitiennes $n\times n$ d{\'e}finies
positives de d{\'e}terminant $1$ et $SU(n)$ le groupe sp{\'e}cial
unitaire d'ordre $n$; consid{\'e}rons la d{\'e}composition polaire
du groupe de Lie r{\'e}el $SL(n, \mathbb{C})$ des matrices
complexes $n\times n$ de d{\'e}terminant $1$,
$$SL(n, \mathbb{C}) = SU(n)\times SH(n).$$
Alors $SH(n)$ muni de la multiplication d{\'e}finie par la
projection sur $SH(n)$ du produit matriciel dans $SL(n,
\mathbb{C})$, est une boucle de Lie mono-alternative {\`a} droite.
Explicitement, la multiplication sur $SH(n)$ est d{\'e}finie par
$$m(a, b) = (ba^{2}b)^{\frac{1}{2}}, \hspace{1cm} \forall a, b \in
SH(n).$$ En effet tout {\'e}l{\'e}ment $d \in SL(n, \mathbb{C})$
se d{\'e}compose de mani{\`e}re unique sous la forme $d = ga$,
o{\`u} $g \in SU(n)$ et $a \in SH(n)$; comme $(\bar{g})^{t}g =
I_{n}$ et $(\bar{a})^{t} = a$ par d{\'e}finition de $SU(n)$ et
$SH(n)$, multipliant $d$ {\`a} gauche  par la conjugu{\'e}e de sa
transpos{\'e}e, on obtient $(\bar{d})^{t}d = a^{2};$ d'o{\`u} $a =
[(\bar{d})^{t}d]^{\frac{1}{2}}$ est la projection de $d$ sur
$SH(n)$. Pour tous $a, b \in SH(n)$, distincts de $I_{n}$ et
$b\neq a^{-1}$, le produit $ab$ n'est {\'e}l{\'e}ment ni de
$SH(n)$, ni de $SU(n)$, c'est un {\'e}l{\'e}ment de $SL(n,
\mathbb{C})$. Comme $(\bar{a})^{t}= a$ et $(\bar{b})^{t}= b$ par
d{\'e}finition de $SH(n)$, de ce qui pr{\'e}c{\`e}de, la
projection du produit $ab$ sur $SH(n)$ est
$$m(a, b) = (ba^{2}b)^{\frac{1}{2}}, \hspace{1cm} \forall a, b \in
SH(n).$$ On v{\'e}rifie sans peine que $(SH(n), m)$ est une boucle
de Lie mono-alternative {\`a} droite.
\end{ex}
De m{\^e}me, l'ensemble des matrices hermitiennes d{\'e}finies
positives $n\times n$ et l'ensemble des matrices sym{\'e}triques
r{\'e}elles d{\'e}finies positives $n\times n$ sont des boucles de
Lie mono-alternatives {\`a} droite, la multiplication $m$
{\'e}tant d{\'e}finie comme dans l'exemple ci-dessus.
\section{Groupes de Lie quasi-doubles}
Les groupes de Lie quasi-doubles sont des g{\'e}n{\'e}ralisations
naturelles des groupes de Lie doubles \cite{Lu-W}.
\begin{defi} Un groupe de Lie quasi-double est un triplet $(G,
G_1, G_2)$, o{\`u}\\
1) G est un groupe de Lie et $G_1$ un sous-groupe de Lie ferm{\'e}
de $G$; \\
2) $G_2$ est une sous-vari{\'e}t{\'e} ferm{\'e}e de $G$ munie de
deux applications analytiques
\begin{center} $\alpha : G_2\times G_2 \longrightarrow G_1$
\hspace{1cm} et \hspace{1cm} $m : G_2\times G_2 \longrightarrow
G_2$ \end{center} satisfaisant l'identit{\'e}
\begin{center} $(ab)_G = \alpha(a, b)m(a, b)$ \hspace{1cm}
$\forall a, b \in G_2$, \end{center} telles que $m$ d{\'e}finit
sur $G_2$ une structure de boucle de Lie mono-alternative {\`a}
droite;\\
3) l'application
\begin{center}
$\begin{array}{ccccc}
\omega &: &G_1\times G_2 &\longrightarrow &G, \\
&  &  (g, a) &\longrightarrow & ga
\end{array}$
\end{center}
est un diff{\'e}omorphisme.\\
\end{defi}
Si $G_1$, $G_2$ et $\omega$ sont d{\'e}finis localement, on dit
que $(G, G_1, G_2)$ est un groupe de Lie quasi-double local.\\
{\bf Remarques :} Si $\alpha(a, b) = e$, $\forall a, b \in G_2,$
o{\`u} $e$ d{\'e}signe l'{\'e}l{\'e}ment neutre de $G$, alors
$G_2$ est aussi un sous-groupe de Lie de $G$ et le triplet $(G,
G_1, G_2)$ d{\'e}finit un groupe de Lie double au sens de
\cite{Lu-W}. Par ailleurs, en identifiant $G$ avec $G_1\times G_2$
par $\omega$, on voit que $\alpha$ et $m$ sont respectivement les
projections sur $G_1$ et $G_2$ de la
restriction de la multiplication de $G$ {\`a} $G_2$. \\
Nous allons {\`a} pr{\'e}sent montrer que si $(G, G_1, G_2)$ est
un groupe de Lie quasi-double, alors il existe deux applications
analytiques
\begin{center} $\chi : G_2\times G_1 \longrightarrow G_1$
\hspace{1cm} et \hspace{1cm}$\sigma :G_2\times G_1 \longrightarrow
G_2$\end{center} telles que $\sigma$  d{\'e}finisse une action du
groupe de Lie $G_1$ sur la vari{\'e}t{\'e} $G_2$ et $G_2$ agisse
sur $G_1$ par $\chi$ modulo l'application $\alpha$. En effet si
$g\in G_1$ et $a\in G_2$, alors $ag\in G$; par cons{\'e}quent
d'apr{\`e}s la condition (3) de la d{\'e}finition
pr{\'e}c{\'e}dente, il existe un unique {\'e}l{\'e}ment $g'\in
G_1$ et un unique {\'e}l{\'e}ment $a'\in G_2$ tels que
 $$ag = g'a'$$
 Posons $a' = a^{g}$ et $g' = g^{a}$, et consid{\'e}rons les
 applications suivantes:
\begin{center}
 $\begin{array}{ccccc}
 \chi & : & G_2\times G_1 &\longrightarrow & G_1,\\
      &   & (a, g) & \longrightarrow & g^{a}\\
 \end{array}$
 \hspace{1cm}
 $\begin{array}{ccccc}
  \sigma & : & G_2\times G_1 & \longrightarrow & G_2\\
         &  & (a, g) & \longrightarrow & a^{g}
 \end{array}$
 \end{center}
 On a le r{\'e}sultat suivant:
 \begin{thm} Soit $(G, G_1, G_2)$ un groupe de Lie quasi-double;
 soient $g, h \in G_1$ et $a,b,c \in G_2$. Alors on a les
 identit{\'e}s suivantes:
 \begin{enumerate}
 \item $(gh)^{a} = g^{a}h^{a^{g}}$
\item $(a^{g})^{h} = a^{gh}$ \item $(g^{b})^{a}\alpha(a^{g^{b}},
b^{g})= \alpha(a, b)g^{m(a, b)}$ \item $[m(a, b)]^{g} =
m(a^{g^{b}}, b^{g})$ \item $\alpha(a, b)\alpha(m(a, b), c) =
[\alpha(b,c)]^{a}\alpha(a^{\alpha(b, c)}, m(b, c))$ \item $m(m(a,
b), c) = m(a^{\alpha(b, c)}, m(b, c))$
\end{enumerate}
Inversement si $G_1$ est un groupe de Lie agissant {\`a} droite
par $\sigma$ sur une boucle de Lie mono-alternative {\`a} droite
$G_2$ dont la loi est d{\'e}finie par $m : G_2\times G_2
\longrightarrow G_2$ et s'il existe deux applications analytiques
$\alpha : G_2\times G_2 \longrightarrow G_1$ et $\chi : G_2\times
G_1 \longrightarrow G_1$ tels que les conditions (1)-(6) soient
v{\'e}rifi{\'e}es, alors il existe une structure de groupe de Lie
sur $G_1\times G_2$ telle que $(G_1\times G_2, G_1, G_2)$
d{\'e}finisse un groupe de Lie quasi-double.
 \end{thm}
 {\bf D{\'e}monstration :} Pour la premi{\`e}re partie, la
 d{\'e}monstration des identit{\'e}s (1)-(6) est une
 cons{\'e}quence directe de l'associativit{\'e} et de
 l'unicit{\'e} de la d{\'e}composition du produit dans $G$.\\
 Pour la seconde partie, on d{\'e}finit la structure de groupe de
 Lie sur $G_1\times G_2$ par :
 $$(g, a)(h, b) = ( gh^{a}\alpha(a^{h}, b), m(a^{h}, b)), \forall g, h \in G_1, \forall a, b \in G_2.$$\\
 Il suffit de v{\'e}rifier l'associativit{\'e} du produit ainsi
 d{\'e}fini sur $G_1\times G_2$, qui est une cons{\'e}quence des
 conditions (1)-(6). $\Box$ \\
{\bf Remarque :} La condition (2) signifie que $\sigma$
 d{\'e}finit une action {\`a} droite du groupe de Lie $G_1$ sur la
 vari{\'e}t{\'e} $G_2$, tandis que la condition (3) signifie
 que $G_2$ agit {\`a} gauche par $\chi$ sur $G_1$ mais que $\chi$
 n'est pas une vraie action {\`a} gauche, le d{\'e}faut d'action {\'e}tant
 mesur{\'e} par l'application $\alpha$; dans la suite,
 nous appellerons une telle application $\chi$ une {\bf action $\alpha-$twist{\'e}e}.
 Si $\alpha(a, b) = e, \forall a, b \in G_2$, alors $G_2$ est un groupe de Lie et $\chi$
devient une action {\`a} gauche de $G_2$ sur la vari{\'e}t{\'e}
$G_1$.
\begin{cor} Soit $(G, G_{1}, G_{2})$ un groupe de Lie
quasi-double. Alors $\forall a, b, c \in G_{2}$, on a :
\begin{enumerate}
\item $\lambda_{m(a, b)}(c) = (\lambda_{\sigma(a, \alpha(b,
c))}o\lambda_{b})(c),$ \item $\rho_{m(a, b)}(c) =
(\rho_{b}o\rho_{a})(\sigma(c, \alpha(a, b)^{-1})),$ \item
$(\rho_{b}o\lambda_{a})(c) = (\lambda_{\sigma(a, \alpha(c,
 b))}o\rho_{b})(c).$
\end{enumerate}
\end{cor}
Ces diff{\'e}rentes relations d{\'e}coulent de la relation (6) du
 th{\'e}or{\`e}me pr{\'e}c{\'e}dent. Si $\alpha(a, b) = e, \forall a, b \in G_2$,
 alors $G_2$ est un groupe de Lie et on retrouve les identit{\'e}s traditionnelles
 sur les tranlations {\`a} gauche et {\`a} droite dans un groupe de Lie.\\
 \begin{ex} Le triplet $(SL(n, \mathbb{C}), SU(n), SH(n))$, o{\`u} $SH(n)$
 est l'ensemble des matrices hermitiennes $n\times n$ d{\'e}finies positives
 de d{\'e}terminant $1$, est un groupe de Lie quasi-double. Les
 applications $\alpha$ et $m$ sur $SH(n)$ sont d{\'e}finies
 respectivement par
 \begin{center}
 $\alpha(a, b) = (ab)(ba^{2}b)^{-\frac{1}{2}}$ \hspace{0.5cm}
 et \hspace{0.5cm} $m(a, b) = (ba^{2}b)^{\frac{1}{2}} \hspace{0.5cm} \forall a, b \in SH(n).$
 \end{center}
 \end{ex}
 D{\'e}finissons {\`a} pr{\'e}sent la notion correspondante pour
 les alg{\`e}bres de Lie en {\'e}tablissant une relation entre les
 deux notions.
\begin{defi} Une alg{\`e}bre de Lie quasi-double  est un triplet
$({\bf g, g_1, g_2})$ o{\`u} ${\bf g}$ est une alg{\`e}bre de Lie,
${\bf g_1}$ une sous-alg{\`e}bre de Lie de ${\bf g}$ et ${\bf
g_2}$ un sous-espace vectoriel de ${\bf g}$ muni de deux
op{\'e}rations bilin{\'e}aires antisym{\'e}triques
\begin{center} $\psi : \Lambda^{2}{\bf g_2} \longrightarrow {\bf g_1}$ \hspace{1cm}
et \hspace{1cm} $\mu : \Lambda^{2}{\bf g_2}\longrightarrow {\bf
g_2}$
\end{center}
tels que
\begin{enumerate}
\item ${\bf g} = {\bf g_1} \oplus {\bf g_2}$ comme somme directe
de sous-espaces vectoriels; \item $[x, y]_{{\bf g}} = \psi(x, y) +
\mu(x, y)$, $\forall x, y \in {\bf g_2}$.
\end{enumerate}
\end{defi}
En un mot,  une alg{\`e}bre de Lie quasi-double est juste une
alg{\`e}bre de Lie qui se d{\'e}compose en somme directe d'une
sous-alg{\`e}bre de Lie et d'un sous-espace vectoriel.\\
{\bf Remarque :} Si $\psi(x, y) = 0, \forall x, y \in {\bf g_2}$,
alors ${\bf g_2}$ est aussi une sous-alg{\`e}bre de Lie de ${\bf
g}$ et le triplet $({\bf g, g_1, g_2})$ d{\'e}finit une
alg{\`e}bre de Lie double au sens de \cite{Lu-W}. Dans
\cite{Al-KS}, les alg{\`e}bres de Lie quasi-doubles $({\bf g, g_1,
g_2})$ telles que ${\bf g}$ est munie d'une forme bilin{\'e}aire
invariante non d{\'e}g{\'e}n{\'e}r{\'e}e par rapport {\`a}
laquelle ${\bf g_1}$ et ${\bf g_2}$ sont des sous-espaces
isotropiques maximaux de ${\bf g}$, sont appel{\'e}es des
quasi-triples de Manin.
\begin{ex} Le triplet $(sl(n, \mathbb{C}), su(n), sh(n))$, o{\`u} $sl(n,
\mathbb{C})$ est l'alg{\`e}bre de Lie des matrices complexes
$n\times n$ de trace nulle, $su(n)$ l'alg{\`e}bre de Lie des
matrices anti-hermitiennes $n\times n$ de trace nulle et $sh(n)$
l'espace vectoriel des matrices hermitiennes $n\times n$ de trace
nulle, est une alg{\`e}bre de Lie quasi-double. Ici, $\psi$ est la
restriction {\`a} $sh(n)$ du crochet de Lie de $sl(n,
\mathbb{C})$, qui est {\`a} valeurs dans $su(n)$ et $\mu = 0$.
\end{ex}
Plus g{\'e}n{\'e}ralement, la d{\'e}composition de Cartan
\cite{He} de toute alg{\`e}bre de Lie semi-simple complexe de
dimension finie, fournit une structure d'alg{\`e}bre de Lie
quasi-double.\\
Le r{\'e}sultat suivant {\'e}tablit une correspondance entre les
groupes de Lie quasi-doubles locaux et les alg{\`e}bres de Lie
quasi-doubles.
\begin{thm} Soit $(G, G_1, G_2)$ un groupe de Lie quasi-double
local; soient ${\bf g, g_1}$ les alg{\`e}bres de Lie de $G$ et
$G_1$ respectivement, et soit ${\bf g_2} = T_{\varepsilon}G_2$
l'espace tangent {\`a} $G_2$ en l'identit{\'e}. Alors $({\bf g,
g_1, g_2})$ est une alg{\`e}bre de Lie quasi-double.\\
Inversement, {\`a} toute alg{\`e}bre de Lie quasi-double, il
correspond  un unique groupe de Lie quasi-double local ({\`a}
isomorphisme pr{\`e}s), dont l'alg{\`e}bre de Lie quasi-double est
celle de d{\'e}part.
\end{thm}
{\bf D{\'e}monstration :} Soit $(G, G_1, G_2)$ un groupe de Lie
quasi-double local; soient ${\bf g, g_1}$ les alg{\`e}bres de Lie
de $G$ et $G_1$ respectivement, et soit ${\bf g_2} =
T_{\varepsilon}G_2$; il est clair que ${\bf g} = {\bf g_1} \oplus
{\bf g_2}$. Pour la conclusion, il suffit de prendre pour $\psi$
le crochet commutateur associ{\'e} {\`a} $\alpha$ et pour $\mu$ le
crochet commutateur associ{\'e} {\`a} $m$.\\
Pour la seconde partie, consid{\'e}rons une alg{\`e}bre de Lie
quasi-double $({\bf g, g_1, g_2})$ et soient $G$, $G_1$ les
groupes de Lie connexes et simplement connexes correspondant {\`a}
${\bf g}$ et ${\bf g_1}$ respectivement. Soit $\Pi : G
\longrightarrow G_1\setminus G$ la projection naturelle de $G$
dans l'espace $G_1\setminus G$ des classes {\`a} droite et $\phi$
la restriction {\`a} ${\bf g_2}$ de l'application compos{\'e}e
$\Pi o (exp) : {\bf g} \longrightarrow G_1\setminus G$. Alors
$\phi$ est un diff{\'e}omorphisme d'un voisinage $V$ de $0$ dans
${\bf g_2}$ dans un voisinage de la classe $\Pi(e)$ dans
$G_1\setminus G$ \cite{He}. Posons $G_2 = expV$ et introduisons
sur $G_2$ une loi de composition interne d{\'e}finie par :
$$m(a, b) = \Pi_{G_2}(ab)$$
$\forall a, b \in G_2$, o{\`u} $ab$ est le produit dans $G$ des
{\'el{\'e}ments $a$ et $b$ et $\Pi_{G_2}$ est la projection locale
sur $G_2$ parall{\`e}lement {\`a} $G_1$, i.e, dans la
d{\'e}composition locale $d=ga \in G$, $g\in G_1$ et $a\in G_2$.
Alors $(G_2, m)$ est une boucle de Lie locale mono-alternative
{\`a} droite \cite{Sa-Mi2}. En effet, par d{\'e}finition de $m$,
on a
$$m(m(a, b^{k}), b^{l}) = \Pi_{G_2}(\Pi_{G_2}(ab^{k})b^{l})=
\Pi_{G_2}((ab^{k})b^{l}) = \Pi_{G_2}(ab^{k+l})$$
et comme pour $b$ suffisamment proche de $\varepsilon$ dans $G_2$
on a $b^{k+l}\in G_2$, alors
$$m(a, m(b^{k}, b^{l})) = \Pi_{G_2}(a\Pi_{G_2}(b^{k}b^{l}) =
\Pi_{G_2}(ab^{k+l})= m(a, b^{k+l}).$$
Par cons{\'e}quent
$$m(m(a, b^{k}), b^{l}) = m(a, b^{k+l})$$
Ainsi $(G_2, m)$ est une boucle de Lie locale mono-alternative
{\`a} droite. D'apr{\`e}s \cite{Sa-Mi2}, cette boucle de Lie locale
mono-alternative {\`a} droite est unique {\`a} isomorphisme pr{\`e}s.\\
Posons $\alpha(a, b) = \Pi_{G_1}(ab)$; de la d{\'e}composition
locale des {\'e}l{\'e}ments de $G$, il r{\'e}sulte que
l'application $ \omega : G_1\times G_2 \longrightarrow G$
d{\'e}finie par $(g, a)\longrightarrow ga$ est un
diff{\'e}omorphisme local et que $(ab)_G = \alpha(a, b)m(a, b)$,
$\forall a, b \in G_2$. En d{\'e}finitive $(G, G_1, G_2)$ est un
groupe de Lie quasi-double local. D'o{\`u} le r{\'e}sultat. $\Box$

Comme pour les groupes de Lie quasi-doubles, nous allons montrer
que pour toute alg{\`e}bre de Lie quasi-double donn{\'e}e $({\bf
g, g_1, g_2})$, il existe deux applications correspondant {\`a}
$\chi$ et $\sigma$. En effet si $\xi \in {\bf g_1}$ et $x \in {\bf
g_2}$, alors $[x, \xi]_{\bf g}\in {\bf g}$; par cons{\'e}quent il
existe $\xi'\in {\bf g_1}$ et $x'\in {\bf g_2}$ tels que
$$[x, \xi] = \xi' + x'.$$
Posons $\xi' = \xi^{x}$ et $x' = x^{\xi}$, et consid{\'e}rons les
applications suivantes :\\
\begin{center}
 $\begin{array}{ccc}
{\bf g_2}\times {\bf g_1}& \longrightarrow & {\bf g_1},\\
(x, \xi)& \longrightarrow & \xi^{x}
\end{array}$
\hspace{1cm}$\begin{array}{ccc}
{\bf g_2}\times {\bf g_1}& \longrightarrow & {\bf g_2}\\
(x, \xi)& \longrightarrow & x^{\xi}
\end{array}$
\end{center}
 On a le r{\'e}sultat suivant
\begin{thm} Soit $({\bf g, g_1, g_2})$ une alg{\`e}bre de Lie
quasi-double. Soient $\xi, \eta, \in {\bf g_1}$, $x, y, z \in {\bf
g_2}$. Alors on a les identit{\'e}s suivantes :
\begin{enumerate}
\item $[\xi, \eta]^{x} = [\xi^{x}, \eta] + [\xi, \eta^{x}] -
(\xi)^{x^{\eta}} + (\eta)^{x^{\xi}}$ \item $x^{[\xi, \eta]} =
(x^{\xi})^{\eta} - (x^{\eta})^{\xi}$ \item $\mu(x, y)^{\xi} =
\mu(x^{\xi}, y) + \mu(x, y^{\xi}) + (x)^{{\xi}^{y}} -
(y)^{\xi^{x}}$ \item $\xi^{\mu(x, y)} = (\xi^{y})^{x} -
(\xi^{x})^{y} + [\xi, \psi(x, y)] + \psi(x^{\xi}, y) + \psi(x,
y^{\xi})$
 \item
$\sum_{\bigcirc}\mu(\mu(x, y), z) = \sum_{\bigcirc}x^{\psi(y, z)}$
\item $\sum_{\bigcirc}\psi(\mu(x, y), z) = \sum_{\bigcirc}\psi(x,
y)^{z}$
\end{enumerate}
Inversement si ${\bf g_1}$ est une alg{\`e}bre de Lie agissant sur
un espace vectoriel ${\bf g_2}$ par $(x, \xi) \longrightarrow
x^{\xi}$, ${\bf g_2}$ muni de deux applications bilin{\'e}aires
antisym{\'e}triques $\mu : \Lambda^{2}{\bf g_2} \longrightarrow
{\bf g_2}$ et $\psi : \Lambda^{2}{\bf g_2} \longrightarrow {\bf
g_1}$, et s'il existe une application de ${\bf g_2}\times {\bf
g_1}$ dans ${\bf g_1}$ d{\'e}finie par $(x, \xi) \longrightarrow
\xi^{x}$, tels que les conditions (1)-(6) ci-dessus soient vraies,
alors il existe une structure d'alg{\`e}bre de Lie sur ${\bf g_1}
\oplus {\bf g_2}$ telle que $({\bf g_1} \oplus {\bf g_2}, {\bf
g_1, g_2})$ d{\'e}finisse une alg{\`e}bre de Lie quasi-double.
\end{thm}
{\bf D{\'e}monstration :} La d{\'e}monstration de la premi{\`e}re
partie du th{\'e}or{\`e}me est une cons{\'e}quence de
l'identit{\'e} de Jacobi dans ${\bf g_1} \oplus {\bf g_2}$.\\
 Pour la seconde partie, on d{\'e}finit
un crochet $[.,.]$ sur ${\bf g_1} \oplus {\bf g_2}$ par :
 $$[\xi, \eta] =  [\xi, \eta]_{\bf g_1}, \forall \xi, \eta \in {\bf
g_1} $$
$$ [x, \xi] =  \xi^{x} +  x^{\xi},  \forall \xi \in {\bf g_1},
\forall x \in {\bf g_2}$$
$$ [x, y]  =  \psi(x, y) + \mu(x, y), \forall x, y \in {\bf g_2}$$
 et on v{\'e}rifit  l'identit{\'e} de
Jacobi pour le crochet ainsi d{\'e}fini, qui est une
cons{\'e}quence des conditions (1)-(6) de la premi{\`e}re partie
du th{\'e}or{\`e}me. Ce qui ach{\`e}ve la d{\'e}monstration.
 $\Box$
\begin{cor} Soit $({\bf g, g_1, g_2})$ une alg{\`e}bre de Lie
quasi-double. Posons $$ [x, y]_{{\bf g_2}} = \mu(x, y),
\hspace{1cm} <x, y, z> = \frac{1}{2}x^{\psi(x, y)}.$$ Alors le
triplet $({\bf g_2}, [.,.]_{{\bf g_2}}, <.,.,.>)$ est une
alg{\`e}bre d'Akivis.
\end{cor}
La d{\'e}monstration du corollaire est une interpr{\'e}tation de
la condition (5) du th{\'e}or{\`e}me pr{\'e}c{\'e}dent.

 Si $(G, G_1, G_2)$ est un groupe de Lie quasi-double, localement, en
appliquant la formule de Campbell-Hausdorff dans $G$ \cite{G}, les
applications $m, \alpha, \sigma$ et $\chi$ admettent les
d{\'e}veloppements limit{\'e}s (d'ordre $3$) suivants au voisinage
de z{\'e}ro :
$$m(x, y) = x + y + \frac{1}{2}\mu(x, y) + \frac{1}{12}\left(\mu(x, \mu(x, y))
+ \mu(y, \mu(y,x)) + x^{\psi(x, y)} + y^{\psi(y, x)}\right)+...$$
$$\alpha(x, y) = \frac{1}{2}\psi(x, y) + \frac{1}{12}\left(\psi(x, \mu(x, y))
+ \psi(y, \mu(y,x)) + \psi(x, y)^{x} + \psi(y, x)^{y}\right)+...$$
$$\sigma(x, \xi) = x + \frac{1}{2}x^{\xi} + \frac{1}{12}\left(\mu(x,
x^{\xi}) + (x)^{\xi^{x}} + (x^{\xi})^{\xi}\right) + ...$$
$$\chi(x, \xi) = \xi + \frac{1}{2}\xi^{x} + \frac{1}{12}\left(\psi(x,
x^{\xi}) + [\xi^{x}, \xi] + (\xi^{x})^{x} + (\xi)^{x^{\xi}}\right)
+...$$ $\forall \xi \in {\bf g_1}, \forall x, y \in {\bf g_2}.$

Soit $(G, G_1, G_2)$ un groupe de Lie quasi-double; pour tout $x
\in T_{\varepsilon}G_{2}$, soit $x^{L}$ le champ de vecteurs
invariant {\`a} gauche sur $G$ associ{\'e} {\`a} $x$,
c'est-{\`a}-dire
$$x^{L}(d) = (T_{e}L_{d})(x)\in T_{d}G, \forall d \in G$$
 et $x^{\lambda}$ l'image directe par $\Pi_{G_2}$ du champ
 $x^{L}$, c'est-{\`a}-dire
$$x^{\lambda}(\Pi_{G_2}(d)) = [T_{\varepsilon}(\Pi_{G_2}\circ
L_{d})](x)\in T_{\Pi({d})}G_2, \forall d \in G.$$
Il est clair que
$x^{\lambda}$ n'est pas invariant {\`a} gauche par les
translations dans $G_2$ du fait de la non associativit{\'e} de
$m$; appelons-le le champ de vecteurs
translat\'e} {\`a} gauche sur $G_2$ associ{\'e} {\`a} $x$.\\
On a le r{\'e}sultat suivant
\begin{prop} Soit $(G, G_1, G_2)$ un groupe de Lie quasi-double; soit
$\rho_{G_2} : T_{e}G_1 \longrightarrow A^{1}(G_2)$,
l'homomorphisme d'alg{\`e}bres de Lie associ{\'e} {\`a} l'action
$\sigma$ de $G_1$ sur $G_2$, $A^{1}(G_2)$ {\'e}tant l'alg{\`e}bre
de Lie des champs de vecteurs sur $G_2$. Alors pour tous $x, y \in
T_{\varepsilon}G_2$ et pour tous $\xi \in T_{e}G$, on a
$$x^{\lambda}(m(a, b)) = (\lambda_{a})_{*}(x^{\lambda}(b)) +
(\rho_{b})_{*}(\rho_{B}(x^{b})(a)),$$
$$\rho_{B}(\xi)(m(a, b))  = (\lambda_{a})_{*}(\rho_{B}(\xi)(b)) +
(\rho_{b})_{*}(\rho_{B}(\xi^{b})(a)),$$
$$[x^{\lambda}, y^{\lambda}] = \rho_{G_2}(\psi(x, y)) + (\mu(x,
y))^{\lambda}$$
 o{\`u} $x^{b}$ et $\xi^{b}$ sont des {\'e}l{\'e}ments de $T_{e}G$,
 d{\'e}finis respectivement par
 \begin{center}
 $x^{b} = \frac{d}{dt}|_{t=0}(\alpha(b, exptx))$ \hspace{1cm} et
 \hspace{1cm} $\xi^{b} = \frac{d}{dt}|_{t=0}(\chi(b, expt\xi))$
 \end{center}
\end{prop}
{\bf D{\'e}monstration :} Ici, on notera $\sigma_{a}(g) =
\sigma(a, g)$. D{\'e}montrons la premi{\`e}re relation; en effet
$\forall a, b \in G_{2}$ et $\forall x \in T_{\varepsilon}G_{2}$, on a:\\
\begin{center}
$\begin{array}{ccc} x^{\lambda}(m(a, b)) & = &
\frac{d}{dt}|_{t=0}(m(m(a, b),
exptx))\\
  &  =  & \frac{d}{dt}|_{t=0}(m(\sigma(a, \alpha(b, exptx)), m(b,
  exptx)))\\
  &  =  & (\Pi_{G_{2}})_{*}(\frac{d}{dt}|_{t=0}(\sigma(a, \alpha(b, exptx)).\Pi_{G_{2}}(b.exptx)))\\
  &  =  & (\Pi_{G_{2}})_{*}((\sigma_{a})_{*}(\frac{d}{dt}|_{t=0}(\alpha(b, exptx))).b)\\
  &     &  + (\Pi_{G_{2}})_{*}((\sigma(a, e).(\frac{d}{dt}|_{t=0}(\Pi_{G_{2}}(b.exptx))))\\
  &  =  & (\Pi_{G_{2}}oR_{b}o\sigma_{a})_{*}(\frac{d}{dt}|_{t=0}(\alpha(b,
  exptx))) + (\Pi_{G_{2}}oL_{a})_{*}(x^{\lambda}(b))\\
  &   =  & (\rho_{b}o\sigma_{a})_{*}(\frac{d}{dt}|_{t=0}(\alpha(b,
  exptx))) + (\lambda_{a})_{*}(x^{\lambda}(b))\\
  &  =  & (\lambda_{a})_{*}(x^{\lambda}(b)) +
  (\rho_{b})_{*}(\rho_{B}(x^{b})(a)).
\end{array}$\\
\end{center}
D'o{\`u} le r{\'e}sultat; en particulier si $\alpha(a, b) = e$,
$\forall a, b \in G_2,$ alors $x^{b} = 0, \forall x \in
T_{\varepsilon}G_{2}$ et on retrouve le fait que $x^{\lambda}$ est
invariant {\`a} gauche sur $G_2$ qui est dans ce cas un groupe de
Lie. La d{\'e}monstration de la deuxi{\`e}me relation est
identique, et celle de la derni{\`e}re d{\'e}coule de
l'identit{\'e} bien connue pour les groupes de Lie en se
pla\c{c}ant dans le groupe de Lie $G$, {\`a} savoir
$$[x^{L}, y^{L}] = [x, y]^{L}.$$
Ce qui ach{\`e}ve la d{\'e}monstration de la proposition. $\Box$
\section{Big{\`e}bres quasi-Lie}
Les big{\`e}bres quasi-Lie sont les objets duaux des
quasi-big{\`e}bres de Lie introduites par Drinfeld \cite{Dr2}.
Plus pr{\'e}cisement
\begin{defi} Une big{\`e}bre quasi-Lie est un quadruplet $({\bf
F}, \mu, \gamma, \psi)$ o{\`u} ${\bf F}$ est un espace vectoriel
de dimension finie sur $\mathbb{K}$, $\mathbb{K} = \mathbb{R}$ ou
$\mathbf{C}$, $\mu : \Lambda^{2}{\bf F}\longrightarrow {\bf F}$,
$\gamma : {\bf F} \longrightarrow \Lambda^{2}{\bf F}$, et $\psi
\in \Lambda^{3}{\bf F^{*}}$ tels que :
\begin{enumerate}
\item $Alt(\gamma\otimes Id)\gamma(x) = 0, \forall x\in {\bf F},$
i.e que $({\bf F}, \gamma)$ est une co-alg{\`e}bre de Lie; \item
$\gamma$ est une d{\'e}rivation par rapport {\`a} $\mu$,
c'est-{\`a}-dire
$$\gamma(\mu(x, y)) = \mu(\gamma(x),
y) + \mu(x, \gamma(y)), \hspace{1cm} \forall x, y \in
T_{\varepsilon}B;$$ \item $Alt(\mu^{t}\otimes Id)\mu^{t}(\xi) =
(\delta_{\gamma}\psi)(\xi), \forall \xi \in {\bf F^{*}}$, o{\`u}
$\mu^{t}$ est l'op{\'e}rateur transpos{\'e} de $\mu$ et
$\delta_{\gamma}$ d{\'e}signe l'op{\'e}rateur cobord d'alg{\`e}bre
de Lie d{\'e}finie par $\gamma$; \item $Alt((\mu^{t}\otimes
Id\otimes Id)\psi) = 0,$
\end{enumerate}
o{\`u} $Alt(x\otimes y\otimes z) = x\otimes y\otimes z + y\otimes
z\otimes x + z\otimes x\otimes y$.
\end{defi}
{\bf Remarques :}\\
a) Dans le cas o{\`u} $\psi = 0$, le triplet $({\bf F}, \mu,
\gamma)$ satisfaisant les conditions ci-dessus est une big{\`e}bre
de Lie \cite{Dr1, KS-Ma, A-KS}. Dans ce cas, la condition (2)
signifie que $\gamma$ est un $1-$cocycle de $({\bf F}, \mu)$ {\`a}
valeurs dans $\Lambda^{2}{\bf F}$ par
rapport {\`a} l'action adjointe d{\'e}finie par $\mu$.\\
b) Si $({\bf F}, \mu, \gamma, \psi)$ est une big{\`e}bre
quasi-Lie, alors $({\bf F^{*}}, \gamma, \mu,  \psi)$ est une
quasi-big{\`e}bre de Lie \cite{KS1, KS2}. Ainsi, comme on le voit,
ces deux notions sont duales l'une de l'autre et ne sont pas
auto-daules.

Soit ${\bf F}$ un espace vectoriel de dimension finie sur
$\mathbb{K}$, $\mathbb{K} = \mathbb{R}$ ou $\mathbf{C}$, $\mu :
\Lambda^{2}{\bf F}\longrightarrow {\bf F}$, $\gamma : {\bf F}
\longrightarrow \Lambda^{2}{\bf F}$, et $\psi \in \Lambda^{3}{\bf
F^{*}}$; en identifiant les applications $\mu$ et $\gamma$ {\`a}
des {\'e}l{\'e}ments de l'alg{\`e}bre ext{\'e}rieure de ${\bf
F^{*}}\oplus {\bf F}$, {\`a} savoir $\mu\in \Lambda^{2}{\bf
F^{*}}\otimes {\bf F}$, $\gamma \in {\bf F^{*}}
\otimes\Lambda^{2}{\bf F}$, et en utilisant les propri{\'e}t{\'e}s
du grand crochet \cite{KS2} nous avons le r{\'e}sultat suivant
\begin{prop} $({\bf F}, \mu, \gamma, \psi)$ est une big{\`e}bre
quasi-Lie si et seulement si $\mu + \gamma + \psi$ est de
carr{\'e} nul par rapport au grand crochet; plus pr{\'e}cisement
si les conditions suivantes sont satisfaites
$$[\gamma, \gamma] = 0$$
$$[\gamma, \mu] =0$$
$$\frac{1}{2}[\mu, \mu] + [\gamma, \psi] = 0$$
$$[\mu, \psi] = 0.$$
\end{prop}
Ce formalisme est plus technique et a l'avantage de simplifier les
calculs dans la th{\'e}orie des big{\`e}bres de Lie et leurs
g{\'e}n{\'e}ralisations \cite{KS2}.
\begin{ex} Soit $({\bf F}, \gamma)$ une co-alg{\`e}bre de Lie;
soit $\omega \in \Lambda^{2}{\bf F^{*}}$ et posons
$$\mu = \delta_{\gamma}(\omega)$$
$$\psi = - \frac{1}{2}[\omega, \omega]^{\gamma}$$
o{\`u} $[.,.]^{\gamma}$ d{\'e}signe le crochet de Schouten
alg{\'e}brique \cite{KS2} associ{\'e} {\`a} $\gamma$. Alors $({\bf
F}, \mu, \gamma, \psi)$ est une big{\`e}bre de quasi-Lie.
\end{ex}
\begin{ex} Consid{\'e}rons la d{\'e}composition polaire de
l'alg{\`e}bre de Lie r{\'e}elle $sl(n, \mathbb{C})$,
c'est-{\`a}-dire
$$sl(n, \mathbb{C}) = su(n) \oplus sh(n).$$
Soient $\gamma$ le crochet de Lie sur $su(n)$ et $\psi$ la
restriction {\`a} $sh(n)$ du crochet de Lie de $sl(n,
\mathbb{C})$, qui est {\`a} valeurs dans $su(n)$. Comme $su(n)$
s'identifie au dual $sh(n)^{*}$ de $sh(n)$ par l'interm{\'e}diaire
de la forme bilin{\'e}aire invariante non
d{\'e}g{\'e}n{\'e}r{\'e}e d{\'e}finie sur $sl(n, \mathbb{C})$ par
$$<x, y> = Im(trace(xy)), \hspace{1cm} x, y \in sl(n, \mathbb{C}),$$
on peut identifier $\psi$ {\`a} un {\'e}l{\'e}ment de
$\Lambda^{3}su(n) \cong \Lambda^{3}(sh(n))^{*}$. Par
cons{\'e}quent $(sh(n), 0, \gamma, \psi)$ est une big{\`e}bre
quasi-Lie.
\end{ex}
Nous allons {\`a} pr{\'e}sent montrer qu'{\`a} toute structure de
big{\`e}bre quasi-Lie sur un espace vectoriel donn{\'e} ${\bf F}$
de dimension finie, il correspond une certaine structure
d'alg{\`e}bre de Lie sur ${\bf F}\oplus {\bf F^{*}}$ et
inversement \cite{KS2, Ba-KS}; faisant ainsi du triplet $({\bf
F}\oplus {\bf F^{*}}, {\bf F^{*}}, {\bf F})$ une alg{\`e}bre de
Lie quasi-double.

Soit $({\bf F}, \mu, \gamma, \psi)$ une big{\`e}bre quasi-Lie;
pour plus de simplicit{\'e}, d{\'e}signons $\gamma$ et sa
transpos{\'e}e par $\gamma$ et identifions $\psi \in
\Lambda^{3}{\bf F^{*}}$ {\`a} une application bilin{\'e}aire
antisym{\'e}trique not{\'e}e aussi par $\psi : \Lambda^{2} {\bf F}
\longrightarrow {\bf F^{*}}$ d{\'e}finie par
$$<\psi(x, y), z> = \psi(x, y, z), \forall x, y, z \in {\bf F}$$
o{\`u} $<.,.>$ d{\'e}signe le produit de dualit{\'e} entre ${\bf
F}$ et ${\bf F^{*}}$. Nous avons le r{\'e}sultat suivant
(\cite{KS2}):
\begin{thm} Le crochet $M$ sur ${\bf F}\oplus {\bf F^{*}}$
d{\'e}fini par :
$$M(x, y) = \mu(x, y) + \psi(x, y), \forall x, y \in {\bf F}$$
$$M(x, \xi) = - ad_{\xi}^{*\gamma}x + ad_{x}^{*\mu}\xi, \forall x
\in {\bf F}, \forall \xi \in {\bf F^{*}}$$
$$M(\xi, \eta) = \gamma(\xi, \eta), \forall \xi, \eta \in {\bf
F^{*}}$$ o{\`u} $ad_{x}^{\mu}y = \mu(x, y)$, $ad_{x}^{*\mu} = -
(ad_{x}^{\mu})^{t}$, $ad_{\xi}^{\gamma}\eta = \gamma(\xi, \eta)$
et $ad_{\xi}^{*\gamma} = - (ad_{\xi}^{\gamma})^{t}$,\\
est un crochet d'alg{\`e}bre de Lie laissant invariant le produit
scalaire canonique d{\'e}fini  sur ${\bf F}\oplus {\bf F^{*}}$ par
$$<x + \xi, y + \eta> = <\xi, y> + <\eta, x>, \forall x, y \in
{\bf F}, \forall \xi, \eta \in {\bf F^{*}}$$
\end{thm}
De fa\c{c}on plus pr{\'e}cise, on montre dans \cite{KS2}, que les
structures de big{\`e}bre quasi-Lie sur un espace vectoriel ${\bf
F}$ sont en correspondance bijective  avec les structures
d'alg{\`e}bre de Lie quasi-double sur $({\bf F}\oplus {\bf F^{*}},
{\bf F^{*}}, {\bf F})$ laissant invariant le produit scalaire
canonique, c'est-{\`a}-dire les structures de quasi-triple de
Manin \cite{Al-KS}. Dans \cite{KS2, Ba-KS}, le couple $({\bf
F}\oplus {\bf F^{*}}, M)$ est appel{\'e} le double de la
big{\`e}bre quasi-Lie $({\bf F}, \mu, \gamma, \psi)$. Ainsi,
d'apr{\`e}s le corollaire 2, {\`a} toute structure de big{\`e}bre
quasi-Lie, il correspond une structure d'alg{\`e}bre d'Akivis.\\
Pour plus de d{\'e}tails sur les big{\`e}bres quasi-Lie, voir
\cite {KS1, KS2, Ba-KS}.
\section{Quasi-groupes de Lie quasi-Poisson}
Dans cette section, nous d{\'e}finissons l'objet
g{\'e}om{\'e}trique correspondant {\`a} une big{\`e}bre quasi-Lie,
puis nous {\'e}tablissons un r{\'e}sultat liant les deux notions.
\begin{defi} Un quasi-groupe de Lie quasi-Poisson est un ensemble
$(B, m, \- G,  \sigma, P, \alpha, \chi, <.,.>)$ o{\`u} $(B, m)$
est une boucle de Lie mono-alternative {\`a} droite, $G$ un groupe
de Lie agissant {\`a} droite sur $B$ par $\sigma : B\times G
\longrightarrow B$, P un champ de bivecteurs s'annulant en
l'identit{\'e} de $B$, $\alpha : B\times B \longrightarrow G$ une
application analytique, $\chi : B\times G \longrightarrow G$ une
action {\`a} gauche $\alpha-$twist{\'e}e  de $B$ sur $G$ et $<.,.>
 : T_{e}G \times T_{\varepsilon}B  \longrightarrow \mathbb{R}$ une
forme bilin{\'e}aire  invariante non d{\'e}g{\'e}n{\'e}r{\'e}e
tels que
\begin{enumerate}
\item $m(m(a, b), c) = m(\sigma(a, \alpha(b, c)), m(b, c)),
\forall a, b, c \in B;$ \item $\alpha(a, b)\alpha(m(a, b), c) =
\chi(a, \alpha(b, c))\alpha(\sigma(a, \alpha(b, c)), m(b, c)),
\forall a, b, c \in B;$ \item $\frac{1}{2}[P, P]_S =
 - (\Lambda^{3}{\rho}_B)(\psi)$, o{\`u} $[.,.]_S$ d{\'e}signe
le crochet de Schouten des champs de multivecteurs \cite{Ko, KS2},
$\rho_B : T_{e}G \longrightarrow A^{1}(B)$ l'homomorphisme
d'alg{\`e}bres de Lie associ{\'e} {\`a} l'action $\sigma$ et
$\psi$ le crochet commutateur associ{\'e} {\`a} $\alpha$; \item
$L_{x^{\lambda}}P = [(L_{x^{\lambda}}P)(\varepsilon)]^{\lambda} -
(\Lambda^{2}\rho_{B})(\psi(x)), \forall x \in T_{\varepsilon}B$,
o{\`u} $x^{\lambda}$ d{\'e}signe le champ de vecteurs
translat{\'e} {\`a} gauche correspondant {\`a} $x$.
\end{enumerate}
\end{defi}
{\bf Remarques :}\\
1) Si $\alpha(a, b) = e, \forall a, b \in B$ (donc $\psi = 0$),
alors $(B, m)$ est un groupe de Lie muni d'une structure de
Poisson; si en plus $B$ est connexe, la condition (4) de la
d{\'e}finition ci-dessus ($\psi = 0$) est {\'e}quivalente {\`a} la
multiplicativit{\'e} \cite{Lu-W, Lu, KS2} de $P$ car
$P(\varepsilon) = 0$, et donc $(B, m, P)$ est un groupe de Lie-Poisson. \\
2) L'invariance de $<.,.>$ signifie que $\forall \xi, \eta \in
T_{e}G$, $\forall x, y \in T_{\varepsilon}B$, on a :\\
a) $<\xi, \mu(x, y)> = <\xi^{y}, x>$, $\mu$ {\'e}tant le commutateur associ{\'e} {\`a} $m$;\\
b) $<[\xi, \eta], x> = <\eta, x^{\xi}>$ o{\`u} $[.,.]$ est le crochet de Lie sur $T_{e}G$; \\
c) $<\psi(x,y), z> = <\psi(y, z), x>.$\\
3) La non d{\'e}g{\'e}n{\'e}rescence  de $<.,.>$ identifie
$T_{e}G$ avec le dual $(T_{\varepsilon}B)^{*}$ de
$(T_{\varepsilon}B)$; d'o{\`u}, de la remarque 2), $\psi \in
\Lambda^{3}(T_{\varepsilon}B)^{*} \cong \Lambda^{3}T_{e}G;$ \\
4) Par d{\'e}finition d'un quasi-groupe de Lie quasi-Poisson, $P$
n'est pas en g{\'e}n{\'e}ral Poisson, mais il d{\'e}finit un
crochet de Poisson sur l'espace $C^{\infty}(B)^{G}$ des fonctions
$G-$invariantes d{\'e}finies sur $B$ \cite{Al-KS}.

On a le r{\'e}sultat suivant
\begin{prop} Soit $(B, m, G, \sigma, P, \alpha, \chi, <.,.>)$
un quasi-groupe de Lie quasi-Poisson (local); posons $\gamma(x) =
(L_{x^{\lambda}}P)(\varepsilon)$, $\forall x \in
T_{\varepsilon}B$. Alors l'application $\gamma : T_{\varepsilon}B
\longrightarrow \Lambda^{2}(T_{\varepsilon}B)$ satisfait les
propri{\'e}t{\'e}s suivantes :
\begin{enumerate}
\item $\gamma^{2} = 0;$ \item $\gamma(\mu(x, y)) = \mu(\gamma(x),
y) + \mu(x, \gamma(y)), \hspace{1cm} \forall x, y \in
T_{\varepsilon}B,$
\end{enumerate}
o{\`u} $\mu$ est le crochet associ{\'e}e {\`a} $m$, et
$\gamma^{2}(x)= (\gamma \otimes Id + Id \otimes
\gamma)(\gamma(x)).$
\end{prop}
{\bf D{\'e}monstration :} 1) En utilisant l'identit{\'e} de Jacobi
gradu{\'e}e du crochet de Schouten, on obtient\\
\begin{center}
$\begin{array}{ccc}
[x^{\lambda}, [P, P]] & = & 2[[x^{\lambda}, P], P] \\
 & = & 2[\gamma(x)^{\lambda} - (\Lambda^{2}\rho_{B})(\psi(x)), P] \\
  & = & 2([\gamma^{2}(x)]^{\lambda} -
  (\Lambda^{3}\rho_{B})(\psi(\gamma(x)))\\
  &   & - [(\Lambda^{2}\rho_{B})(\psi(x)), P]), \forall x\in
T_{\varepsilon}B.
\end{array}$
\end{center}
 Par ailleurs, de la condition (3)
de la d{\'e}finition d'un quasi-groupe de Lie quasi-Poisson, on a
\begin{center}
$\begin{array}{ccc} [x^{\lambda}, [P, P]](\varepsilon) & = & -
2[x^{\lambda}, \Lambda^{3}\rho_{B}(\psi)](\varepsilon) =
- 2(L_{x^{\lambda}}\Lambda^{3}\rho_{B}(\psi))(\varepsilon)\\
&  = &  - 2\frac{d}{dt}|_{t=0}((\Lambda^{3}\rho_{B}(\psi)(exptx))exp(-tx))\\
 &  = &  - 2\frac{d}{dt}|_{t=0}((\Lambda^{3}(\rho_{exp(-tx)}o\sigma_{exptx}))_{*}(\psi))
\end{array}$
\end{center}
o{\`u} $\sigma_{a}(g) = \sigma(a, g), \forall a\in B, \forall g
\in G$. Comme $\psi \in \Lambda^{3}(T_{\varepsilon}B)^{*}$ et
$(\sigma_{\varepsilon})_{*} =0$, la r{\`e}gle de d{\'e}rivation de
la composition des applications prouve que $[x^{\lambda}, [P,
P]](\varepsilon) = 0.$\\
Ainsi, en {\'e}valuant l'{\'e}galit{}\'e
$$[x^{\lambda}, [P, P]] =
2([\gamma^{2}(x)]^{\lambda} -
(\Lambda^{3}\rho_{B})(\psi(\gamma(x))) -
[(\Lambda^{2}\rho_{B})(\psi(x)), P])$$ en l'identit{\'e}, sachant
que $P(\varepsilon) = 0$ et $\rho_{B}$ s'annule en l'identit{\'e},
on trouve $\gamma^{2}(x) = 0,
\forall x\in T_{\varepsilon}B$; d'o{\`u} $\gamma^{2} = 0.$\\
2) En utilisant la relation (4) de la d{\'e}finition d'un
quasi-groupe de Lie quasi-Poisson et la relation $[x^{\lambda},
y^{\lambda}] = \rho_{B}(\psi(x, y) + (\mu(x, y))^{\lambda}$, on
obtient
\begin{center}
$\begin{array}{ccc} \gamma(\mu(x, y)) & = & (L_{\mu(x,
y)^{\lambda}}P)(\varepsilon)\\
& = & (L_{[x^{\lambda}, y^{\lambda}]}P)(\varepsilon) - (L_{\rho_{B}(\psi(x,y))}P)(\varepsilon) \\
 & = &  (L_{x^{\lambda}}L_{y^{\lambda}}P)(\varepsilon) -
 (L_{y^{\lambda}}L_{x^{\lambda}}P)(\varepsilon)\\
 & = & [x^{\lambda}, [y^{\lambda}, P]](\varepsilon) - [y^{\lambda}, [x^{\lambda}, P]](\varepsilon)\\
 & = & [x^{\lambda}, \gamma(y)^{\lambda}](\varepsilon) - [x^{\lambda}, (\Lambda^{2}\rho_{B})(\psi(y))](\varepsilon)\\
 &   & - [y^{\lambda}, \gamma(x)^{\lambda}](\varepsilon) + [y^{\lambda}, (\Lambda^{2}\rho_{B})(\psi(x))](\varepsilon)\\
 & = & \mu(\gamma(x), y) + \mu(x, \gamma(y)).
 \end{array}$
 \end{center}
D'o{\`u}
$$\gamma(\mu(x, y)) = \mu(\gamma(x),
y) + \mu(x, \gamma(y)), \hspace{1cm} \forall x, y \in
T_{\varepsilon}B.$$ Ce
qui ach{\`e}ve la d{\'e}monstration de la proposition. $\Box$\\
{\bf Remarques :} La relation $\gamma^{2} = 0$ signifie que  la
lin{\'e}arisation  de $P$ \cite{W} en l'identit{\'e} de $B$ est un
co-crochet de Lie sur $T_{\varepsilon}B$ ou encore un crochet de
Lie sur $(T_{\varepsilon}B)^{*}\cong T_{e}G$ qu'on va identifier
avec le crochet de Lie de $T_{e}G$. Ainsi l'invariance de $<.,.>$
prouve que
\begin{center} $\xi^{x} = ad_{x}^{*\mu}\xi, \hspace{1cm}
x^{\xi} = - ad_{\xi}^{*\gamma}x$
\end{center}
La relation (2) de la proposition exprime le fait que $\gamma$ est
une d{\'e}rivation par rapport {\`a} $\mu$. En termes de grand
crochet, les deux relations s'{\'e}crivent respectivement comme
suit
$$[\gamma, \gamma] = 0$$
$$[\gamma, \mu] = [\mu, \gamma] = 0.$$
Le r{\'e}sultat suivant {\'e}tablit une correspondance entre les
quasi-groupes de Lie quasi-Poisson locaux et les big{\`e}bres
quasi-Lie.
\begin{thm} L'espace tangent en l'identit{\'e} {\`a} tout
quasi-groupe de Lie quasi-Poisson local $(B, m, G, \- \sigma, P,
\alpha, \chi, <.,.>)$ est muni d'une structure de big{\`e}bre
quasi-Lie $({\bf F}, \mu, \gamma, \psi)$, o{\`u} ${\bf F} =
T_{\varepsilon}B$, $\mu$ le crochet commutateur d{\'e}fini {\`a}
partir de la loi de composition $m$, $\gamma = d_{\varepsilon}P$
la lin{\'e}arisation de $P$ en l'identit{\'e} et $\psi$ le crochet
commutateur associ{\'e} {\`a} l'application $\alpha$.

Inversement, {\`a} toute big{\`e}bre quasi-Lie r{\'e}elle $({\bf
F}, \mu, \gamma, \psi)$, il correspond  un unique quasi-groupe de
Lie quasi-Poisson local ({\`a} isomorphisme pr{\`e}s) dont la
big{\`e}bre quasi-Lie tangente co{\"\i}ncide avec la big{\`e}bre
quasi-Lie donn{\'e}e.
\end{thm}
{\bf D{\'e}monstration :} Soit $(B, m, G, \sigma, P, \alpha, \chi,
<.,.>)$ un quasi-groupe de Lie quasi-Poisson local. Soit ${\bf F}
= T_{\varepsilon}B$; soient $\mu$ et $\psi$ les crochets
commutateurs  associ{\'e}s {\`a} $m$ et $\alpha$ respectivement.
Comme $P(\varepsilon)= 0$, on peut consid{\'e}rer la
lin{\'e}arisation de $P$ en $\varepsilon$ et d{\'e}signons-la par
$\gamma$. D'apr{\`e}s la proposition 3, $\gamma$ est un co-crochet
de Lie sur ${\bf F}$ et une d{\'e}rivation par rapport {\`a} $\mu$.\\
Consid{\'e}rons {\`a} pr{\'e}sent la relation $$m(m(a, b), c) =
m(\sigma(a, \alpha(b, c)), m(b, c)), \forall a, b, c \in B.$$
Alors en effectuant localement un d{\'e}veloppement limit{\'e}
d'ordre 3 au voisinage de z{\'e}ro des deux membres de cette
{\'e}galit{\'e}, on obtient pour tous $x, y, z \in {\bf F}$
\begin{center}
$\begin{array}{c} 3\mu(\mu(x, y), z) + \mu(x, \mu(y, z)) +
\mu(y, \mu(x, z)) + y^{\psi(x, z)} +
x^{\psi(y, z)} - 3z^{\psi(x, y)} = \\
= 3\mu(x, \mu(y, z)) + \mu(y,\mu(z, x)) + \mu(z, \mu(y, x)) +
y^{\psi(z, x)} + 3x^{\psi(y, z)} + z^{\psi(y, x)}
\end{array}$
\end{center}
En utilisant l'antisym{\'e}trie des applications $\mu$ et $\psi$,
on obtient
$$\sum_{\bigcirc}\mu(\mu(x, y), z) = \sum_{\bigcirc} x^{\psi(y, z)}.$$
Mais l'invariance de $<.,.>$ prouve que
$$x^{\psi(y, z)}= - ad_{\psi(y, z)}^{*\gamma}x;$$
d'o{\`u}
$$\sum_{\bigcirc}\mu(\mu(x, y), z) = - \sum_{\bigcirc}ad_{\psi(y, z)}^{*\gamma}x, \hspace{1cm} \forall x, y, z \in {\bf F}.$$
Par transposition, on trouve que cette {\'e}galit{\'e} est
{\'e}quivalente {\`a}
$$Alt(\mu^{t}\otimes Id)\mu^{t}(\xi) = (\delta_{\gamma}\psi)(\xi), \hspace{1cm} \forall \xi \in {\bf F^{*}}.$$
Enfin consid{\'e}rons la relation
$$\alpha(a, b)\alpha(m(a, b), c) =
\chi(a, \alpha(b, c))\alpha(\sigma(a, \alpha(b, c)), m(b, c)),
\forall a, b, c \in B;$$ comme pr{\'e}c{\'e}demment, en effectuant
localement un d{\'e}veloppement limit{\'e} d'ordre 3 au voisinage
de z{\'e}ro des deux membres de cette {\'e}galit{\'e}, on obtient
 $$\sum_{\bigcirc}\psi(\mu(x, y), z) = \sum_{\bigcirc}\psi(y, z)^{x}, \hspace{1cm} \forall x, y, z \in {\bf F}.$$
 L'invariance de $<.,.>$ prouve que
 $$\psi(y, z)^{x} = ad_{x}^{*\mu}\psi(y, z);$$
 d'o{\`u}
 $$\sum_{\bigcirc}\psi(\mu(x, y), z) = \sum_{\bigcirc}ad_{x}^{*\mu}\psi(y, z), \hspace{1cm} \forall x, y, z \in {\bf F}.$$
 Par transposition, on trouve que cette derni{\`e}re
 {\'e}galit{\'e} est {\'e}quivalente {\`a}
$$Alt((\mu^{t}\otimes Id\otimes Id)\psi) = 0.$$
En d{\'e}finitive, $({\bf F}, \mu, \gamma, \psi)$ est une
big{\`e}bre quasi-Lie.

Inversement soit $({\bf F}, \mu, \gamma, \psi)$ une big{\`e}bre
quasi-Lie; d'apr{\`e}s le th{\'e}or{\`e}me 4, il lui correspond
une structure d'alg{\`e}bre de Lie quasi-double sur $({\bf
F}\oplus {\bf F^{*}}, {\bf F^{*}}, {\bf F})$ laissant invariant le
produit scalaire canonique. Soit $(D, G, B)$ le groupe de Lie
quasi-double local correspondant {\`a} $({\bf F}\oplus {\bf
F^{*}}, {\bf F^{*}}, {\bf F})$ comme d{\'e}crit par le
th{\'e}or{\`e}me 2. $\forall a, b \in B$, $\forall g \in G$ posons
$$m(a, b) = \Pi_{B}(ab), \alpha(a, b) =
\Pi_{G}(ab), \sigma(a, g) = \Pi_{B}(ag), \chi(a, g) = \Pi_{G}(ag)
.$$ Il est clair que $\sigma$ est une action de $G$ sur $B$ et
$\chi$ est une action $\alpha$-twist{\'e}e. De la relation
d'associativit{\'e} dans $D$, $(a(bc)_{D})_{D} = ((ab)_{D}c)_{D}$,
$\forall a, b, c \in B$, on obtient les relations (1) et (2) de la
d{\'e}finition d'un quasi-groupe de Lie quasi-Poisson.\\
Par ailleurs, d'apr{\`e}s \cite{KS2}, par la donn{\'e}e de la
big{\`e}bre de quasi-Lie $({\bf F}, \mu, \gamma, \psi)$, ${\bf
F}\oplus {\bf F^{*}}$ est munie d'une structure de
quasi-big{\`e}bre de Lie quasitriangulaire $({\bf F}\oplus {\bf
F^{*}}, M, \Gamma, \psi)$, o{\`u} $M$ est la structure
d'alg{\`e}bre de Lie sur ${\bf F}\oplus {\bf F^{*}}$ d{\'e}crite
par le th{\'e}or{\`e}me 4, $\Gamma$ est un co-crochet sur ${\bf
F}\oplus {\bf F^{*}}$ d{\'e}fini par $\Gamma = \gamma - \mu^{t} -
\psi$; explicitement
$$\Gamma(x, \xi) = \gamma(x) - \mu^{t}(\xi) - \psi(x) \in
\Lambda^{2}({\bf F}\oplus {\bf F^{*}})$$ o{\`u} $\gamma : {\bf F}
\longrightarrow \Lambda^{2}{\bf F}$, $\mu^{t} : {\bf F^{*}}
\longrightarrow \Lambda^{2}{\bf F^{*}}$ (transposition), $\psi :
{\bf F}\longrightarrow \Lambda^{2}{\bf F^{*}}$.\\
De la th{\'e}orie des quasi-big{\`e}bres de Lie et groupes de Lie
quasi-Poisson \cite{KS2}, on sait que cette structure sur ${\bf
F}\oplus {\bf F^{*}}$ s'int{\`e}gre en une structure de groupe de
Lie quasi-Poisson sur $D$ not{\'e}e par
$(D, P, \psi)$, o{\`u}\\
$\bullet$ $P$ est un champ de bivecteurs multiplicatif \cite{Lu-W,
Lu, KS2} sur le groupe de Lie $D$, c'est-{\`a}-dire
$P(gh) = (L_{g})_{*}P(h) + (R_{h})_{*}P(g), \forall g, h \in D$; \\
$\bullet$ $\frac{1}{2}[P, P]_{S} = \psi^{R} - \psi^{L}$, o{\`u}
$\psi^{L}$ (respectivement $\psi^{R}$) d{\'e}signe le champ de
trivecteurs
invariant {\`a} gauche (respectivement {\`a} droite) sur $D$ associ{\'e} {\`a} $\psi$;\\
$\bullet$ $[P, \psi^{L}]_{S} = [P, \psi^{R}]_{S} = 0.$\\
Consid{\'e}rons {\`a} pr{\'e}sent l'image directe par $\Pi_{B} : D
\longrightarrow B$ du champ de bivecteurs $P$ et notons-la par
$P_{B}$, c'est-{\`a}-dire
$$P_{B}(\Pi_{B}(d)) = (\Lambda^{2}T_{d}\Pi_{B})(P(d)),
\hspace{2cm} \forall d \in D.$$ $P_{B}$ est bien un champ de
bivecteurs d{\'e}fini sur $B$ avec $P_{B}(\varepsilon) = 0$ o{\`u}
$\varepsilon = \Pi_{B}(e)$, {\`a} cause du caract{\`e}re
multiplicatif de $P$. Par ailleurs, prenant l'image directe par
$\Pi_{B}$ des deux membres de l'{\'e}galit{\'e}
$$\frac{1}{2}[P, P]_{S} = \psi^{R} - \psi^{L},$$
on obtient
$$\frac{1}{2}[P_{B}, P_{B}]_{S} = (\Lambda^{3}(\Pi_{B})_{*})(\psi^{R}) - (\Lambda^{3}(\Pi_{B})_{*})(\psi^{L}).$$
Mais $(\Lambda^{3}(\Pi_{B})_{*})(\psi^{R}) = 0$ car $\psi \in
\Lambda^{3}(T_{\varepsilon}B)^{*}$ et $\Pi_{B}\circ R_a =
\rho_{a}\circ \Pi_{B}$ par construction de la structure du groupe
de Lie quasi-double $(D, G, B)$; d'autre part
$$(\Lambda^{3}(\Pi_{B})_{*})(\psi^{L}) =
(\Lambda^{3}\rho_{B})(\psi).$$ D'o{\`u}
$$\frac{1}{2}[P_{B}, P_{B}]_{S} = - (\Lambda^{3}\rho_{B})(\psi).$$
Enfin, comme $D$ est connexe, la multiplicativit{\'e} de  $P$ est
{\'e}quivalente {\`a} la relation
$$L_{X^{L}}P = [(L_{X^{L}}P)(e)]^{L}, \forall X \in T_{e}D = {\bf F}\oplus {\bf
F^{*}},$$ o{\`u} $X^{L}$ d{\'e}signe le champ de vecteurs
invariant {\`a} gauche  sur $D$ associ{\'e} {\`a} $X$. Mais $
(L_{X^{L}}P)(e) = \Gamma(X) = (\gamma - \mu^{t}
- \psi)(X), \forall X \in T_{e}D = {\bf F}\oplus {\bf F^{*}}.$\\
Consid{\'e}rons la relation ci-dessus avec $X = x \in {\bf F}$,
c'est-{\`a}-dire
$$L_{x^{L}}P = [\Gamma(x)]^{L}= [\gamma(x)- \psi(x)]^{L};$$
comme $\psi(x) \in \Lambda^{2}{\bf F^{*}}$, en prenant l'image
directe par $\Pi_{B}$ des deux membres de cette {\'e}galit{\'e},
on obtient
$$L_{x^{\lambda}}P_{B} = [\gamma(x)]^{\lambda} - (\Lambda^{2}\rho_{B})(\psi(x)) =
[(L_{x^{\lambda}}P_{B})(\varepsilon)]^{\lambda} -
(\Lambda^{2}\rho_{B})(\psi(x)).$$ Pour  la forme bilin{\'e}aire
invariante non d{\'e}g{\'e}n{\'e}r{\'e}e, il suffit de prendre
pour $<.,.>$ la restriction {\`a} ${\bf F}\times {\bf F^{*}}$ du
produit scalaire canonique d{\'e}fini sur ${\bf F}\oplus {\bf F^{*}}$.\\
En d{\'e}finitive, $(B, m, G, \sigma, P_{B}, \alpha, \chi, <.,.>)$
est un quasi-groupe de Lie quasi-Poisson local dont la big{\`e}bre
quasi-Lie tangente est bien celle de d{\'e}part. Ce qui ach{\`e}ve
la d{\'e}monstration du th{\'e}or{\`e}me. $\Box$
\begin{cor} Soit $({\bf F}, \mu, \gamma, \psi)$ une big{\`e}bre
de quasi-Lie et soit $(B, m, G, \- \sigma, P_{B}, \alpha, \chi,
<.,.>)$ le quasi-groupe de Lie quasi-Poisson local correspondant
comme d{\'e}crit par le th{\'e}or{\`e}me pr{\'e}c{\'e}dent. Alors
le champ de bivecteurs $P_{B}$ sur $B$ satisfait la relation
$$L_{\rho_{B}(\xi)}P_{B} = - (\Lambda^{2}\rho_{B})(\mu^{t}(\xi)),
\hspace{1cm} \forall \xi \in T_{e}G.$$
\end{cor}
{\bf D{\'e}monstration :} En effet, par d{\'e}finition
$$P_{B}(\Pi_{B}(d)) = (\Lambda^{2}T_{d}\Pi_{B})(P(d)),
\hspace{2cm} \forall d \in D.$$ D'apr{\`e}s le th{\'e}or{\`e}me
pr{\'e}c{\'e}dent, $P$ satisfait la relation
$$L_{X^{L}}P = [\Gamma(X)]^{L}, \forall X \in T_{e}D = {\bf F}\oplus {\bf
F^{*}},$$ o{\`u} $X^{L}$ d{\'e}signe le champ de vecteurs
invariant {\`a} gauche sur $D$ associ{\'e} {\`a} $X$ et $\Gamma(X)
= (\gamma - \mu^{t} - \psi)(X)$.\\
Pour $X = \xi \in T_{e}G = (T_{\varepsilon}B)^{*}$, prenant
l'image directe par $\Pi_{B}$ des deux membres de la relation
pr{\'e}c{\'e}dente sachant que $(\Pi_{B})_{*}(\xi^{L}) =
\rho_{B}(\xi)$, $\mu^{t}(\xi) \in \Lambda^{2}T_{e}G$, $\gamma :
{\bf F} \longrightarrow \Lambda^{2}{\bf F}$, et $\psi : {\bf F}
\longrightarrow \Lambda^{2}{\bf F^{*}}$, on obtient la relation
anonc{\'e}e.$\Box$\\
{\bf Remarques :} \\
a) Si $\alpha(a, b) = e,$ $\forall a, b \in B$ et si le groupe de
Lie $G$ est connexe, la relation {\'e}nonc{\'e}e dans le
corollaire pr{\'e}c{\'e}dent traduit de mani{\`e}re
{\'e}quivalente le fait que $\sigma$ est une action de Poisson \cite{Lu-W}.\\
b) Dans \cite{Al-KS, B-C}, le couple $(B, P_{B})$ v{\'e}rifiant la
condition (3) de la d{\'e}finition d'un quasi-groupe de Lie
quasi-Poisson et la relation {\'e}nonc{\'e}e dans le corollaire
pr{\'e}c{\'e}dent, est appel{\'e} un $G-$espace quasi-Poisson,
o{\`u} $G$ est muni de sa structure de groupe de Lie quasi-Poisson
induite par celle de $D$. Dans \cite{Al-KS}, une construction
similaire du champ de bivecteurs $P_{B}$ est donn{\'e}e sur
l'ensemble $D/G$ des classes {\`a} gauche, et l'action de $G$ sur
$D/G$ est une action {\`a} gauche; d'o{\`u} l'apparition de
diff{\'e}rences de signes par rapport {\`a} la n{\^o}tre,
notamment dans l'expression de $P_{B}$ et de son carr{\'e} par
rapport au crochet de Schouten.
\begin{ex} Tout groupe de Lie-Poisson $(G, P)$ est un quasi-groupe de Lie
quasi-Poisson; il suffit de consid{\'e}rer son groupe de
Lie-Poisson dual $G^{*}$ avec lequel ils agissent l'un sur l'autre
\cite{Lu-W} et de prendre $\psi = 0$.
\end{ex}
\begin{ex} Consid{\'e}rons la d{\'e}composition polaire du groupe
de Lie r{\'e}el $SL(2, \mathbb{C})$ des matrices complexes
$2\times 2$ de d{\'e}terminant $1$,
$$SL(2, \mathbb{C}) = SU(2)\times SH(2),$$
o{\`u} $SH(2)$ est l'ensemble des matrices hermitiennes $2\times
2$ d{\'e}finies positives de d{\'e}terminant $1$ et $SU(2)$ le
groupe sp{\'e}cial unitaire d'ordre $2$. Soit
\begin{displaymath}
e_{1}= \frac{1}{\sqrt{2}}\left(\begin{array}{cc}
1 & 0 \\
0 & -1 \\
\end{array}\right),
e_{2}= \frac{1}{\sqrt{2}}\left(\begin{array}{cc}
0 & 1 \\
1 & 0 \\
\end{array}\right),
 e_{3}= \frac{1}{\sqrt{2}}\left(\begin{array}{cc}
0 & i \\
-i & 0 \\
\end{array}\right)
\end{displaymath}
la base canonique de $sh(2)$, l'espace tangent {\`a} $SH(2)$ en
l'identit{\'e},  et
\begin{displaymath}
\varepsilon_{1}= \frac{1}{\sqrt{2}}\left(\begin{array}{cc}
i & 0 \\
0 & -i \\
\end{array}\right),
\varepsilon_{2}= \frac{1}{\sqrt{2}}\left(\begin{array}{cc}
0 & i \\
i & 0 \\
\end{array}\right),
 \varepsilon_{3}= \frac{1}{\sqrt{2}}\left(\begin{array}{cc}
0 & -1 \\
1 & 0 \\
\end{array}\right).
\end{displaymath}
la base duale de $su(2)$ relativement {\`a} la forme
bilin{\'e}aire
$$<x, \xi> = Im(trace(x\xi)), \hspace{1cm} \forall x \in sh(2),
\forall \xi \in su(2).$$ Par rapport au crochet de Lie de $sl(2,
\mathbb{C})$, nous avons les relations de commutation suivantes
dans les deux bases
$$[e_1, e_2] = - \sqrt{2}\varepsilon_3, \hspace{1cm} [e_1, e_3] =
\sqrt{2}\varepsilon_2, \hspace{1cm} [e_2, e_3] = -
\sqrt{2}\varepsilon_1$$
$$[\varepsilon_1, \varepsilon_2] = \sqrt{2}\varepsilon_3, \hspace{1cm} [\varepsilon_1, \varepsilon_3] =
- \sqrt{2}\varepsilon_2, \hspace{1cm} [\varepsilon_2,
\varepsilon_3] = \sqrt{2}\varepsilon_1$$
$$[e_1, \varepsilon_2] = \sqrt{2}e_3, \hspace{1cm} [e_1, \varepsilon_3] =
- \sqrt{2}e_2, \hspace{1cm} [e_2, \varepsilon_1] = -
\sqrt{2}e_3,$$
$$[e_2, \varepsilon_3] =
\sqrt{2}e_1, \hspace{1cm} [e_3, \varepsilon_1] = \sqrt{2}e_2,
\hspace{1cm} [e_3, \varepsilon_2] = - \sqrt{2}e_1,$$
$$[e_i, \varepsilon_i] = 0, \hspace{2cm} \forall i = 1, 2, 3.$$
Pour tous $a, b \in SH(2)$ et tous $g\in SU(2)$, posons
$$m(a, b) = \Pi_{SH(2)}(ab)= (ba^{2}b)^{\frac{1}{2}}, \hspace{1cm} \alpha(a, b) = \Pi_{SU(2)}(ab)= (ab)(ba^{2}b)^{-\frac{1}{2}},$$
$$\sigma(a, g) = \Pi_{SH(2)}(ag) = (g^{-1}a^{2}g)^{\frac{1}{2}}, \hspace{0.5cm} \chi(a, g) = \Pi_{SU(2)}(ag)= (ag)(g^{-1}a^{2}g)^{-\frac{1}{2}}.$$
Consid{\'e}rons le champ de bivecteurs sur $SH(2)$ d{\'e}fini par
$$P = \frac{1}{2}\sum_{i=1}^{3}e_{i}^{\lambda}\wedge
\rho_{SH(2)}(\varepsilon_{i}),$$ o{\`u} $e_{i}^{\lambda}$ est le
champ de vecteurs translat{\'e} {\`a} gauche (par rapport {\`a} la
loi $m$) associ{\'e} {\`a} $e_i$ et $\rho_{SH(2)}$ est
l'homomorphisme d'alg{\`e}bres de Lie associ{\'e} {\`a}
$\sigma$, l'action de $SU(2)$ sur $SH(2)$.\\
Alors $(SH(2), m, SU(2), \sigma, P, \alpha, \chi, <.,.>)$ est un
quasi-groupe de Lie quasi-Poisson dont la big{\`e}bre quasi-Lie
tangente est $(sh(2), 0, \gamma, \psi)$, o{\`u} $\gamma$ est le
crochet de Lie sur $su(2)$ et $\psi$ la restriction du crochet de
$sl(2, \mathbb{C})$ aux {\'e}l{\'e}ments de $sh(2)$. Le groupe de
Lie quasi-Poisson dual est $(SU(2), 0, \psi)$.
\end{ex}
En effet, les conditions (1) et (2) de la d{\'e}finition d'un
quasi-groupe de Lie quasi-Poisson sont {\'e}videntes, elles
s'obtiennent par un calcul direct simple utilisant les
d{\'e}finitions des diff{\'e}rentes applications; il est clair
aussi que $P(\varepsilon) = 0$ par d{\'e}finition de $P$. Par
ailleurs, les calculs utilisant les relations de commutation entre
les {\'e}l{\'e}ments des deux bases, montrent que
$$\frac{1}{2}[P, P]_S = - (\Lambda^{3}\rho_{SH(2)})(\psi) = 0,$$
o{\`u} $\psi$ consid{\'e}r{\'e} comme {\'e}l{\'e}ment de
$\Lambda^{3}su(2) \cong \Lambda^{3}(sh(2))^{*}$ est {\'e}gal {\`a}
$- \sqrt{2}(\varepsilon_1\wedge \varepsilon_2\wedge
\varepsilon_3)$ d'apr{\`e}s les relations de commutation; ce qui
nous donne la condition (3) de la d{\'e}finition. V{\'e}rifions
{\`a} pr{\'e}sent la condition (4); en effet par un calcul direct
utilisant les relations de commutation ci-dessus, on trouve
\begin{center}
$\begin{array}{ccc}
L_{e_{1}^{\lambda}}P & =
&\sqrt{2}(e_2\wedge e_3)^{\lambda} +
\sqrt{2}(\Lambda^{2}\rho_{SH(2)})(\varepsilon_2\wedge \varepsilon_3)\\
& = &[(L_{e_{1}^{\lambda}}P)(\varepsilon)]^{\lambda} -
(\Lambda^{2}\rho_{SH(2)})(\psi(e_1))
\end{array}$
\end{center}
\begin{center}
$\begin{array}{ccc}
L_{e_{2}^{\lambda}}P & = & -
\sqrt{2}(e_1\wedge e_3)^{\lambda} -
\sqrt{2}(\Lambda^{2}\rho_{SH(2)})(\varepsilon_1\wedge
\varepsilon_3) \\
 & = & [(L_{e_{2}^{\lambda}}P)(\varepsilon)]^{\lambda}
- (\Lambda^{2}\rho_{SH(2)})(\psi(e_2))
\end{array}$
\end{center}
\begin{center}
$\begin{array}{ccc}
L_{e_{3}^{\lambda}}P & = & \sqrt{2}(e_1\wedge
e_2)^{\lambda} +
\sqrt{2}(\Lambda^{2}\rho_{SH(2)})(\varepsilon_1\wedge
\varepsilon_2) \\
 & = & [(L_{e_{3}^{\lambda}}P)(\varepsilon)]^{\lambda}
- (\Lambda^{2}\rho_{SH(2)})(\psi(e_3)).
\end{array}$
\end{center}
En d{\'e}finitive, $(SH(2), m, SU(2), \sigma, P, \alpha, \chi,
<.,.>)$ est un quasi-groupe de Lie quasi-Poisson.\\
{\bf Remarque :} Dans cet exemple, $P$ est Poisson mais n'est pas
multiplicatif au sens de \cite{Lu-W, Lu, KS2}, car $\alpha(a, b)$
n'est pas la matrice identit{\'e} d'ordre 2 pour tous $a, b \in
SH(2)$; une telle structure d{\'e}finit un quasi-groupe de
Lie-Poisson. Plus g{\'e}n{\'e}ralement, si
$(\Lambda^{3}\rho_{B})(\psi) = 0$  avec $\psi \neq 0$, on dit que
$(B, m, G, \- \sigma, P_{B}, \alpha, \chi, <.,.>)$ est un
quasi-groupe de Lie-Poisson. La construction de la structure
d{\'e}finie ci-dessus sur $SH(2)$ s'{\'e}tend naturellement {\`a}
$SH(n)$ pour $n \geq 3$.

\vspace{1cm}

{\bf Remerciements :} L'auteur remercie sinc{\`e}rement Madame le
Professeur Y. Kosmann-Schwarzbach pour ses remarques et
suggestions sur le contenu du travail. Les remerciements vont
{\'e}galement {\`a} la Section Math{\'e}matiques de Abdus Salam
International Centre for Theoretical Physics (ICTP, Trieste,
Italie) pour l'hospitalit{\'e} et le soutien financier, sans
lesquels la r{\'e}alisation du pr{\'e}sent document n'aurait pas
{\'e}t{\'e} possible. Enfin, l'auteur remercie le R{\'e}seau
Africain de G{\'e}om{\'e}trie et d'Alg{\`e}bre Appliqu{\'e}es au
D{\'e}veloppement (RAGAAD) pour avoir financer son voyage {\`a}
travers CIMPA/SARIMA.

\end{document}